\documentclass[12pt,reqno]{amsart}
\usepackage{amssymb, amsfonts, amsbsy, latexsym, color}
\usepackage[dvips]{graphicx}

\newcommand{\bitem}{\begin{itemize}}
\newcommand{\eitem}{\end{itemize}}
\newcommand{\beq}{\begin{equation}}
\newcommand{\eeq}{\end{equation}}
\newcommand{\beqn}{\begin{eqnarray*}}
\newcommand{\eeqn}{\end{eqnarray*}}
\newcommand{\goto}{\rightarrow}
\newcommand{\cC}{{\mathcal C}}
\newcommand{\cA}{{\mathcal A}}

\newcommand{\cS}{{\mathcal S}}

\newcommand{\cN}{{\mathcal N}}
\newcommand{\cQ}{{\mathcal Q}}
\newcommand{\bR}{{\bf R}}

\newcommand{\eps}{{\epsilon}}
\newcommand{\argmin}{\mbox{argmin}}

\def\t{\tilde}
\newcommand{\ip}[2]{\langle#1,#2\rangle}

\newcommand{\absip}[2]{| \langle#1,#2\rangle |}
\newcommand{\norm}[1]{\|#1\|}
\def\supp{{\text{\rm supp}}}

\def\cH{\mathcal{H}}

\def\ZZ{\mathbb{Z}}
\def\bZ{\mathbb{Z}}

\def\RR{\mathbb{R}}
\def\bR{\mathbb{R}}

\def\cT{\mathcal{T}}

\def\t{\tilde}
\def\eps{\varepsilon}
\def\epsilon{\varepsilon}

\newtheorem{theorem}{Theorem}[section]
\newtheorem{corollary}{Corollary}[section]
\newtheorem{proposition}{Proposition}[section]
\newtheorem{lemma}{Lemma}[section]
\newtheorem{definition}{Definition}[section]
\newtheorem{example}{Example}[section]

\begin{document}

\title[Clustered Sparsity and Separation of Cartoon and Texture]{Clustered Sparsity and\\ Separation of Cartoon and Texture}

\author[G. Kutyniok]{Gitta Kutyniok}
\address{Department of Mathematics, Technische Universit\"at Berlin,  10623 Berlin, Germany}
\email{kutyniok@math.tu-berlin.de}

\thanks{The author would like to thank David Donoho for various discussions on this and related topics.
She is grateful to the Department of Statistics at Stanford University and the Department
of Mathematics at Yale University for their hospitality and support during her visits.
The author acknowledges support by the Einstein
Foundation Berlin, by Deutsche Forschungsgemeinschaft (DFG) Heisenberg fellowship KU 1446/8,
Grant SPP-1324 KU 1446/13 and DFG Grant KU 1446/14, and by the DFG Research Center {\sc Matheon}
``Mathematics for key technologies'' in Berlin.}

\begin{abstract}
Natural images are typically a composition of cartoon and texture structures.  A medical image
might, for instance, show a mixture of gray matter and the skull cap. One common task is to separate such an image
into two single images, one containing the cartoon part and the other containing the texture part.
Recently, a powerful class of algorithms using sparse approximation and $\ell_1$ minimization
has been introduced to resolve this problem, and numerous inspiring empirical results have already
been obtained.

In this paper we provide the first thorough theoretical study of the separation of a combination of cartoon
and texture structures in a model situation using this class of algorithms. The methodology we
consider expands the image in a combined dictionary consisting of a curvelet tight frame and a
Gabor tight frame and minimizes the $\ell_1$ norm on the analysis side. Sparse approximation properties
then force the cartoon components into the curvelet coefficients and the texture components into the Gabor
coefficients, thereby separating the image. Utilizing the fact that the coefficients are {\em clustered
geometrically}, we prove that at sufficiently fine scales arbitrarily precise separation is
possible. Main ingredients of our analysis are the novel notion of {\em cluster coherence} and
{\em clustered/geometric sparsity}. Our analysis also provides a deep understanding on when separation
is still possible.
\end{abstract}


\maketitle


\section{Introduction}
\label{sec:intro}

Natural images are typically a composition of `cartoon` and `texture`. Think, for instance, of
a medical image in which we might see a mixture of gray matter and the skull cap, or, more general,
a mixture of tissue and bones.
Often, it is essential to very cleanly separate the cartoon from the texture part, i.e., to generate
two single images from the original one, for separate analysis. However, each educated person these
days would say that it is entirely impossible to solve a problem with only one known datum and two unknowns.
Intriguingly, the deep reason for why this separation is possible in the considered situation is the
stark morphological difference between both structures.

The main idea of the empirical results exploiting applied harmonic analysis methodologies (see, e.g.,
\cite{SED03,SED04,SED05,SMBED05,ESQD05} -- in comparison to  PDE-based separation
methods (see exemplary  \cite{VO03} and references therein), use
 the fact that there exist sparsifying dictionaries for cartoons and for textures. Dictionary
learning on natural images indicates that a curvelet (or shearlet) system might be best
adapted to the cartoon part, whereas a Gabor system might be best adapted to the texture
part (cf. \cite{OF96}). Roughly speaking, the image is then expanded into a combined
dictionary of curvelets and a Gabor system, and the $\ell_1$ norm on the analysis side
is minimized. Sparse approximation properties then force the cartoon components into the curvelet
coefficients and the texture components into the Gabor coefficients, thereby separating the image.
However, no theoretical results are available to date which provide a deep mathematical understanding
of why separation is possible.

In the paper \cite{DK08a} (see also \cite{Kut12}), the apparently similar problem of separating point- and curve-like
structures was considered, a distributional model was developed, and clustered/geometric sparsity
and cluster coherence were introduced to derive an asymptotic separation result; asymptotic in
the scale. However, the analysis required for the situation of cartoon and texture differs
significantly, since a Gabor system does not exhibit a scaling component, and a distributional
model is not feasible for the texture part. Our analysis will though still be based on the
novel viewpoint of clustered/geometric sparsity and cluster coherence.

There exists an intriguing connection with quasicrystals \cite{Gou05}, since those objects
can be regarded as both periodic (texture-like) as well as non-periodic with sharp boundaries
(cartoon-like). This raises the question of how periodization of a cartoon can make it
separable from a single cartoon. This highly fascinating question is though beyond the scope
of this paper, but one goal for future exploration.


\subsection{Model of Cartoon Part}

Intuitively, cartoons are smooth image parts separated from other areas by an edge.
The first model of cartoons has been introduced in \cite{CD04}, and this is what we
intend to use also here. The basic idea is to choose a closed boundary curve
and then fill the interior and exterior part with $C^2$ functions. Here, we will
adapt the definition from \cite{CD04} slightly.

\begin{figure}[ht]
\begin{center}
\vspace*{0.3cm}
\includegraphics[height=1.25in]{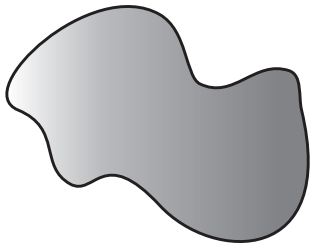}
\vspace*{-0.3cm}
\end{center}
\caption{Example of a cartoon-like image.}
\label{fig:Cartoon}
\end{figure}

We define a model $\cC$ of a cartoon as follows: Let $\tau$ be a closed, non-intersecting, and
regular curve in $\bR^2$, i.e., $\tau$ is $C^\infty$ and of finite arclength, and let
$B_\tau$ denote the interior of $\tau$. Then we choose $\cC$ as
\[
\cC = f_0 + f_1 \cdot 1_{B_\tau},
\]
where $f_0, f_1 \in L^2(\bR^2) \cap C^2(\bR^2)$ with compact support. An example of such a cartoon
is illustrated in Figure \ref{fig:Cartoon}.

To set up the analysis in a reasonable way, we will need to restrict our model of cartoons slightly,
which will be done in Subsection \ref{subsec:balancing}.


\subsection{Model of Texture Part}
\label{subsec:texture}

People have debated for years over an appropriate model for the texture content of
an image. The outcome can be seen nowadays in an extensive variety of texture
models both for the deterministic as well as statistical setting. We would
like to remind the reader of models such as Gaussian random fields,
or the $H^1$-model \cite{OSV03}. All of those emphasize different features of
texture content such as for the $H^{-1}$-model a special oscillatory behavior
`dual' to cartoons, which are sometimes modeled as $BV$-functions, is highlighted.
However, certainly even an intersection of all these models is far from being
capable of accurately describing texture, especially since the perception of texture
also differs from observer to observer. In the area of computer vision very interesting
studies have been performed, for instance by Zucker \cite{BZ03,DZ01a,DZ01b}, which we
would like to bring to the reader's attention.

Our understanding of texture is a function which is concentrated on a small part
of the image with a periodic structure. This justifies the use of a very `classical'
approach by modeling texture as a function which is sparse in a Gabor frame.
However, we would like the reader to be aware of the fact that certainly also this
deterministic model is quite simplified and possible extensions could be incorporated
in it, which we would like to briefly discuss:
\bitem
\item {\em Warping}: Let us consider an image containing a table cloth such as the famous
`Barbara' image. Remembering that this table cloth was in fact folded, we could think of
improving our model by using diffeomorphism to model this warping effect, which was done
in \cite{CM02}.
\item {\em Affine transform}: Again abusing `Barbara' as an example, we recall that the
table cloth was in addition also sheared. This could be handled by application of an
affine transform to the texture content.
\item {\em Windowing}: Our running example also reminds us that texture -- here a pattern
of a table cloth -- could be contained in a polygon, e.g., in particular, the pattern is broken.
Thus refining our model further would force us to consider windowed `Gabor' texture models,
which, in particular, requires a careful treatment of the boundaries. This also illustrates
how delicate the distinction of texture from cartoons can be, and let us wonder where
in such a transition to draw the line between texture and cartoons.
\item {\em Textons}: Finally, textures commonly do not appear as singletons, but might
also overlap. This can be attacked by considering a model of texture which is a linear
combination of shifted `Gabor' texture model with coefficients being iid to decouple the
texture parts. As a reference for this general concept which in literature is normally
referred to as textons, we would like to mention \cite{ZGWX05}.
\eitem

We now define our model for texture as follows: Let $g \in L^2(\bR^2)$ be a window with
$|g(x)| \sim e^{-|x|}$, $\hat{g} \in C^\infty(\bR^2)$, and frequency support $\supp \; \hat{g} \subseteq [-1,1]^2$
satisfying
\[
\sum_{n \in \ZZ} |\hat{g}(\xi + n)|^2 = 1,\qquad \xi \in \bR^2.
\]
For a `size parameter' $s > 0$, let $g_s$ denote the $L^2$-normalized scaled version of $g$ defined by
\[
g_s(x) = s \cdot g(s x).
\]
With this choice, $\hat{g}_s ( \xi) = s^{-1} \cdot \hat{g}(s^{-1} \xi)$,
and hence
\[
\supp \; \hat{g}_s \subseteq [-s,s]^2
\]
as well as
\[
\sum_{n \in \ZZ} |\hat{g}_s(\xi + s n)|^2 = s^{-2},\qquad \xi \in \bR^2.
\]
Further, let $(d_{m,n})_{m,n \in \bZ^2}$ be a sequence of complex numbers belonging
to $\ell_1 \cap \ell_2$. Then the function of interest to us as a model for texture
is defined by
\[
\cT_s = \sum_{m,n \in \bZ^2} d_{m,n} g_s(x-\tfrac{m}{2s}) e^{i (s n)' x}.
\]


\subsection{Model of Image Composed of a Cartoon- and Texture-Part}

Having introduced both constituents, we now assume that the image
\beq \label{eq:trueproblem}
f_s = \cC + \cT_s
\eeq
is observed, however both the cartoon component $\cC$ as well as the texture component
$\cT_s$ are unknown to us. Thus we face the task of extracting those from $f_s$, which
is what we will analyze in this paper.

For the purpose of matching the energy of cartoon- and texture-part so that the separation
becomes non-trivial, in Subsection \ref{subsec:balancing}, we will link $s$ to a particular
other parameter.


\subsection{Sparsifying Overcomplete Systems}
\label{subsec:systems}

Since we intend to utilize $\ell_1$ minimization to separate $\cC$ and $\cT_s$
in the sense of extracting $\cC$ and $\cT_s$ from $f_s$,
appropriate systems need to be selected which provide (relatively) sparse
representations for those components. Of interest to us are the following two
systems:
\bitem
\item {\em Curvelets} -- a highly directional tight frame with increasingly anisotropic
elements at fine scales.
\item {\em Gabor tight frame} -- a tight frame with time-frequency-balanced elements.
\eitem

These systems are constructed as follows. Using the function $g$ introduced in Subsection
\ref{subsec:texture}, we first define the Gabor tight frame at
spatial position index $m$ and frequency position index $n$ by the Fourier transform
\[
(\hat{g}_s)_\lambda (\xi) = \hat{g}_s(\xi- s n) e^{i \frac{m}{2s}' \xi},
\]
where we let $\lambda = (m,n)$ index spatial and frequency position. The bandsize is
indexed by the parameter $s$.
It follows from \cite{DGM86} that this system constitutes a tight frame for $L^2(\bR^2)$.

To define curvelets, we let $W$ be the inverse Fourier transform of a wavelet, where
$W$ belongs to $C^\infty(\bR)$ and is compactly supported on $[1/2,2]$, for instance,
suitably scaled Lemari\`{e}-Meyer wavelets possess these properties. Further, we
choose a ``bump function'' $V \in C^\infty$, which is compactly supported on $[-1,1]$.
We then define {\em continuous curvelets} at scale $a>0$, orientation $\theta \in [0,\pi)$, and spatial position $b \in \bR^2$
by their Fourier transforms
\[
      \hat{\gamma}_{a,b,\theta}(\xi) =  a^{\frac{3}{4}}  \cdot W(a|\xi|) V(a^{-1/2}(\omega-\theta))
         \cdot e^{i b' \xi}.
\]
See \cite{CD04,CD05a} for more details.
The {\em curvelet tight frame} is then (essentially) defined as a sampling of $b$ on a
series of regular lattices
\[
   \{ R_{\theta_{j,\ell}} D_{a_j}  \bZ^2 \}, \qquad j \geq j_0, \quad \ell = 0, \dots, 2^{\lfloor j/2 \rfloor} -1 ,
\]
where $R_{\theta}$ is planar rotation by $\theta$ radians, $a_j = 2^{-j}$,
$\theta_{j,\ell} = \pi \ell / 2^{j/2}$, $\ell = 0, \dots, 2^{j/2}-1$,
and  $D_a$ is anisotropic dilation by $diag(a,\sqrt{a})$, i.e., the curvelets
at scale $j$, orientation $\ell$, and spatial position $k = (k_1,k_2)$ are
given by the Fourier transform
\[
      \hat{\gamma}_{\eta}(\xi) =  2^{-j\frac{3}{4}}  \cdot W(|\xi|/2^{j}) V((\omega-\theta_{j,\ell})2^{j/2})
         \cdot e^{i (R_{\theta_{j,\ell}}D_{2^{-j}}k)' \xi},
\]
where $\eta = (j,k,\ell)$ index scale, orientation, and scale.
(For a precise statement, see \cite[Section 4.3, pp. 210-211]{CD05b}).

Using the {\it same} window $W$, we can construct a family of filters $F_j$ with transfer functions
\[
    \hat{F}_j(\xi) = W(|\xi|/2^{j}), \qquad \xi \in \bR^2 .
\]
These filters allow us to decompose a function $f$ into pieces $f_j$ with different scales,
the piece $f_j$ at subband $j$ arises from filtering $f$ using $F_j$:
\[
   f_j = F_j \star f;
\]
the Fourier transform $\hat{f}_j$ is supported in the annulus with inner radius $2^{j-1}$ and outer radius
$2^{j+1}$. Because of our assumption on $W$, we can reconstruct the original function from these
pieces using the formula
\[
   f = \sum_j F_j \star f_j, \qquad f \in L^2(\bR^2) .
\]

This allows us to split the seeked decomposition of the observed image $f_s$ (see \eqref{eq:trueproblem}) into infinitely many
decompositions depending on scale $j$ by setting
\[
  \cC_j = F_j \star \cC \quad \mbox{and} \quad \cT_{s,j} = F_j \star \cT_s,
\]
and, for each $j$,  considering
\beq \label{eq:scaleproblem}
f_{s,j} = \cC_j + \cT_{s,j}.
\eeq
For each $j$, we assume that $f_{s,j}$ is known to us, and we aim to compute $\cC_j$ and $\cT_{s,j}$ from it.
By the choice of the filter, we can then derive $\cC$ and $\cT_s$ from those.


\subsection{Frequency Matching}
\label{subsec:balancing}

In order to set up the separation problem in a reasonable way, the two filtered components  $\cC_j$ and $\cT_{s,j}$
have to be comparable as we go to  finer scales, so that the ratio of energies is more or less independent
of scale. This ensures that separation is challenging at {\it every} scale.
Thus we need to compute the norms  $\| \cC_{j}\|_2$ and  $\| \cT_{s,j}\|_2$, and then link $s$ to $2^{j}$ in
such a way that both quantities are comparable. This will then also allow us to drive an asymptotic analysis
based on the scale $j$ -- upon which the size $s$ then depends.

Computing the asymptotic behavior of the norm of our filtered cartoon model, we face the problem that the behavior of
$\hat{\cC}$ is not always the same asymptotically for all possible elements. To resolve this
problem, we restrict the model slightly to cartoons which satisfy
\beq \label{eq:extraconditiononC}
|\hat{\cC}(\xi)| \sim |\xi|^{-3/2}, \qquad \xi \in \RR^2.
\eeq
We wish to note that this is only a very mild restriction, since for the boundary of the sphere $\cC = S^2$, it
can be proven that this condition is satisfied (see, e.g., \cite{CD05a}). In general, Randol \cite{Ran69} (cf. also
\cite{Hla50}) proved, for instance, that if $\cC \in C^4$ and if the Gaussian curvature is non-zero in every point
of $\partial \cC$, then $|\hat{\cC}(\xi)| \le C \cdot |\xi|^{-3/2}$, $\xi \in \RR^2$.

With this additional assumption, we can now derive the following lemma whose proof is outsourced to Subsection \ref{subsec:proofs_1}.

\begin{lemma}
 \label{lem:estimateforCj}
\[
\| \cC_j\|_2^2 \sim  2^{-j}, \qquad j \goto \infty.
\]
\end{lemma}

We next compute the asymptotic behavior of the norm of our filtered texture model and also refer to
Subsection \ref{subsec:proofs_1} for the proof.

\begin{lemma}
 \label{lem:estimateforTs}
\[
\| \cT_{s,j}\|_2^2 \sim  \sum_{m, \tilde{m}} e^{-\frac{|m-\tilde{m}|}{2}} \sum_{n \in \ZZ^2 \cap \cA_{s,j}} d_{m,n} \overline{d_{\tilde{m},n}},
\]
where
\[
\cA_{s,j} = \{ \xi \in \RR^2 : \tfrac{2^{j-1}}{s} \le |\xi| \le \tfrac{2^{j+1}}{s}\}.
\]
\end{lemma}

The {\em energy matching condition} on $s$ and $j$ can now be derived from Lemmata \ref{lem:estimateforCj} and \ref{lem:estimateforTs}
by choosing $s = s_j$ to satisfy
\beq \label{eq:energybalancing}
\sum_{m, \tilde{m}} e^{-\frac{|m-\tilde{m}|}{2}} \sum_{n \in \ZZ^2 \cap \cA_{s_j,j}} d_{m,n} \overline{d_{\tilde{m},n}} = 2^{-j}.
\eeq
We now rewrite \eqref{eq:scaleproblem} as
\beq \label{eq:scaleproblem2}
f_{j} = \cC_j + \cT_{j},
\eeq
and also simply write $\cA_j$ instead of $\cA_{s_j,j}$.

For illustrative purposes we present the energy matching condition for one particular sequence $(d_{m,n})_{m,n \in \bZ^2}$.

\begin{example}
\label{exa:specialdmn}
If we exemplarily choose
\[
|d_{m,n}| \sim |m|^{-(2+\delta)} \cdot |n|^{-(2+\delta)}, \qquad \delta > 0,
\]
and observe that for the asymptotics we can ignore the `side terms' $m \neq \tilde{m}$, the energy matching condition \eqref{eq:energybalancing}
becomes
\beq \label{eq:exampleenergymatching}
2^{-j} \sim \sum_m |m|^{-(4+2\delta)} \cdot \sum_{n \in \ZZ^2 \cap \cA_{j}} |n|^{-(4+2\delta)}.
\eeq
Now $\sum_m |m|^{-(4+2\delta)}$ is constant, and since
\[
\int_0^{2\pi} \int_{\frac{2^{j-1}}{s_j}}^{\frac{2^{j+1}}{s_j}} r^{-(4+2\delta)} \cdot r dr d\varphi \sim \left(\frac{2^j}{s_j}\right)^{-2(1+\delta)},
\]
from \eqref{eq:exampleenergymatching}, we finally obtain the condition $2^{-j} \sim \left(2^j/s_j\right)^{-2(1+\delta)}$,
hence we can choose $s_j$ as
\beq \label{eq:energybalancingexample}
s_j = 2^{\frac{1+2\delta}{2+2\delta}j}.
\eeq
Thus, in particular, we have $s_j = \omega(2^{j/2})$ and $s_j = o(2^j)$ in this case. This implies that energy matching is
achieved, if the size of a repeated patch of the texture, i.e., $1/{s_j}$, is `slightly' larger than the thickness
of the filtered boundary of the cartoon, more precisely, between $2^{-j}$ and $2^{-j/2}$. Hence intuitively separation
seems doubtful in this case. It will thus be surprising that in Corollary \ref{coro:main}, we will show arbitrarily
precise separation with the scale $j \to \infty$, if only $\delta > 1$.
\end{example}


\subsection{Separation via $\ell_1$ Minimization}

In our analysis we aim at exploring the geometrical contents and, in particular, their
difference of the cartoonlike structure $\cC$ and texturelike structure $\cT_s$. In
\cite{DK08a}, we have introduced a fundamental notion which allows to derive estimates
on the accuracy of $\ell_1$ based separation of morphological objects. It was therein
applied to prove that, given an image composed of pointlike and curvelike structures
at all sufficiently fine scales, nearly-perfect separation can be achieved. Here we
are concerned with a different situation, since such an asymptotical result is
not possible, simply due to the fact that Gabor systems don't have a scaling parameter.
Although our analysis will therefore be different in nature, the common bracket with
\cite{DK08a} is again the utilization of the geometry of the components. This allows
us to apply some of the novel techniques developed in \cite{DK08a} also in the
setting considered in this paper. For the convenience of the reader, we will briefly
recall those in this subsection.

Suppose we have two tight frames $\Phi_1$, $\Phi_2$ in a Hilbert space $\cH$, and a
signal vector $S = S_1^0 + S_2^0 \in \cH$. We then consider the following optimization problem
\begin{equation} \label{PSep}
    (\mbox{Sep}) \qquad    (S_1^\star, {S}_2^\star)
        = \argmin_{S_1,S_2}  \| \Phi_1^T S_1 \|_1 + \| \Phi_2^T S_2 \|_1
          \mbox{ subject to } S = S_1 + S_2.
\end{equation}
We wish to remark that intentionally the norm is placed on the {\em analysis} coefficients
rather than on the {\em synthesis} coefficients as in Basis Pursuit \cite{CDS01} to avoid
self-terms in the frame setting. Hence we might not recover the most sparsest expansion,
but it will turn out that this expansion is nevertheless sufficient for component separation.
Also it is algorithmically not feasible to optimize over all possible expansions.
In many studies of $\ell_1$ optimization, one then considers the mutual coherence
\[
   \mu (\Phi, \Psi) = \max_{j}  \max_{i} | \langle \phi_{i}, \psi_{j} \rangle|,
\]
whose importance was shown by \cite{DH01}. This might be thought of as singleton coherence.
In contrast, to exploit the knowledge of the geometrical structure of the components -- more
precisely, the location of the `large' coefficients of their frame expansions --, in \cite{DK08a}
the notion of cluster coherence was introduced, which bounds coherence between a single member
of a frame $\Psi$ and a cluster of members of a second frame $\Phi$, clustered at $\cS$.

\begin{definition}
Given tight frames
$\Phi=(\phi_i)_i$ and $\Psi=(\psi_j)_j$
and an index subset $\cS$ associated with
expansions in frame $\Phi$,
we define the \emph{cluster coherence}
\[
        \mu_c (\cS, \Phi; \Psi) = \max_{j}  \sum_{i \in \cS} | \langle \phi_{i}, \psi_{j} \rangle|.
\]
\end{definition}

With this notion at hand, the following estimate concerning the accuracy of
the extraction of $S_1^0$ and $S_2^0$ from the signal $S$ holds true:

\begin{proposition}\cite[Props. 1 + 2]{DK08a} \label{prop:coherenceestimate}
Suppose that $S$ can
be decomposed as $S=S_1^0+S_2^0$ so that
each component $S_i^0$ is relatively sparse in $\Phi_i$, $i=1,2$, i.e.,
\[
\norm{1_{\cS_1^c} \Phi_1^T S_1^0}_1 + \norm{1_{\cS_2^c} \Phi_2^T S_2^0}_1
\le \delta.
\]
Let $(S_1^\star,S_2^\star)$ solve (\ref{PSep}).
Then
\[
  \norm{S_1^\star-S_1^0}_2 + \norm{S_2^\star-S_2^0}_2
\le \frac{2\delta}{1-2\mu_c},
\]
where
\[
\mu_c = \max( \mu_c(\cS_1,\Phi_1;\Phi_2), \mu_c(\cS_2,\Phi_2;\Phi_1)).
\]
\end{proposition}

The concepts of this section will now be applied to  (\mbox{\sc Sep}), at scale $j$ only. For this,
the tight frames are $\Phi_1$, the full curvelet tight frame, and $\Phi_2$, the full Gabor tight frame,
and $S$ is our filtered observed image $f_{j}$, which satisfies (cf. \eqref{eq:scaleproblem2}),
\[
f_{j} =  \cC_j + \cT_{j}.
\]
We apply the optimization problem (\mbox{\sc Sep}), getting subsignal components $S_1^\star$ and
$S_2^\star$, which we then relabel as the curvelet component $C_j$ and Gabor component $T_j$.

The cluster of indices of both, sparse approximation error and cluster coherence, depends on the
choice of the cluster indices $\cS_1$ -- now denoted by $\cS_{1,j}$  -- for the curvelet system
and the cluster indices $\cS_2$ -- now denoted by $\cS_{2,j}$ -- for the Gabor system.
The sparse approximation error $\delta$ and the cluster coherence $\mu_c$ -- let us remind the
reader that these are merely analysis tools and are not part of the $\ell_1$ minimization procedure --
relevant for a scale $j$ will then be denoted by
\[
\delta_{1,j} \quad \mbox{and} \quad \delta_{2,j}
\]
for relative sparsity for the cartoon and texture part, respectively, as well as
\[
(\mu_c)_{1,j} = \mu_c(\cS_{1,j}, \{ \gamma_\eta\}; \{(g_s)_\lambda\}) \quad \mbox{and} \quad (\mu_c)_{2,j} = \mu_c(\cS_{2,j}, \{(g_s)_\lambda\}; \{ \gamma_\eta\})
\]
for the cluster coherence for a cluster of curvelets and Gabor elements, respectively.


\subsection{Asymptotic Separation Result}

Our aim is an asymptotic separation result for the relative $L^2$-error of the purported cartoon part $C_j$
and texture part $T_j$ decays to zero as $j \to \infty$. By Proposition \ref{prop:coherenceestimate}, this follows
if the clusters $\cS_{1,j}$ and $\cS_{2,j}$ can be chosen such that
\[
\max(\delta_{1,j},\delta_{2,j}) = o(2^{-j/2}), \quad \mbox{as } j \to \infty,
\]
and
\[
\max((\mu_c)_{1,j},(\mu_c)_{2,j}) \to 0, \quad \mbox{as } j \to \infty.
\]
Our main result shows that this is indeed possible. For this, the set of significant coefficients $\cS_{2,j}$ for $\cT_{j}$ is defined by
\[
\cS_{2,j} = B(0,r_{1,j}) \times B(0,r_{2,j}),
\]
for some $r_{1,j}, r_{2,j} > 0$, where $B(0,r)$ shall denote the closed $\ell_2$ ball around the origin in $\RR^2$.

\begin{theorem} \label{theo:main}
Assuming energy matching \eqref{eq:energybalancing}, we set
\[
\cA_{j} = \{ \xi \in \RR^2 : \tfrac{2^{j-1}}{s_j} \le |\xi| \le \tfrac{2^{j+1}}{s_j}\}
\]
and
\[
M_{a_j,\theta,s_j} = R_\theta \cdot s_j^{-1} (\text{\rm supp }\hat{\gamma}_{a_j,0,0} + B_1(0,1)) \cap \bZ^2.
\]
Suppose that there exist $r_{1,j}, r_{2,j} > 0$ which satisfy
\[
\sum_{\stackrel{(m,n) \not\in B(0,r_{1,j}) \times B(0,r_{2,j})}{n \in \cA_{j}}} \sum_{\tilde{m} \in \ZZ^2} |d_{\tilde{m},n}| e^{-\frac{|\tilde{m}-m|}{2}}
= o(2^{-j/2}), \qquad j \to \infty
\]
and
\begin{enumerate}
\item[(i)] in case $s_j = o(2^{j})$ as $j \to \infty$,
\[
|B(0,r_{2,j}) \cap M_{a_j,0,s_j} \cap \ZZ^2| = o(2^{3j/4} \cdot s_j^{-1}),\quad \mbox{as } j \to \infty,
\]
\item[(ii)] as well as in case $s_j = \Omega(2^{j})$ as $j \to \infty$,
\[
|B(0,r_{1,j}) \cap \ZZ^2| = o(2^{-3j/4} \cdot s_j),\quad \mbox{as } j \to \infty.
\]
\end{enumerate}
Then we have asymptotically near-perfect separation:
\[
    \frac{ \|C_j - \cC_j \|_2 + \|T_{j} - \cT_{j} \|_2}{\|\cC_{j}\|_2 + \|\cT_{j}\|_2} \goto 0, \qquad j \goto \infty.
\]
\end{theorem}

Most interesting is the consideration of the special case already focussed on in Example \ref{exa:specialdmn}.
In this situation asymptotically near-perfect separation is always attained provided that $\delta > 1$. The
precise statement is the following

\begin{corollary} \label{coro:main}
Suppose that the sequence $(d_{m,n})_{m,n \in \bZ^2}$ satisfies
\[
|d_{m,n}| \sim |m|^{-(2+\delta)} \cdot |n|^{-(2+\delta)}, \qquad \delta > 1.
\]
Then we have asymptotically near-perfect separation:
\[
    \frac{ \|C_j - \cC_j \|_2 + \|T_{j} - \cT_{j} \|_2}{\|\cC_{j}\|_2 + \|\cT_{j}\|_2} \goto 0, \qquad j \goto \infty.
\]
\end{corollary}

Theorem \ref{theo:main} and Corollary \ref{coro:main} are both proved in Subsection \ref{subsec:proofs_1}.


\subsection{Interpretation}
\label{subsec:intuition}

\begin{figure}[ht]
\begin{center}
\vspace*{0.3cm}
\includegraphics[height=1.25in]{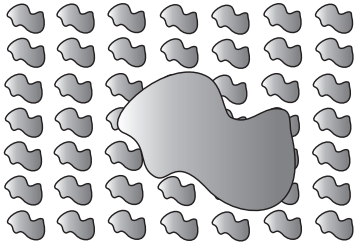}
\vspace*{-0.3cm}
\end{center}
\caption{Periodic small cartoons versus one large cartoon.}
\label{fig:periodiccartoons}
\end{figure}

To heuristically understand why cartoon and texture can be separated at all, consider the worst case scenario
that an image is composed of one cartoon and a periodized cartoon as is illustrated in Figure \ref{fig:periodiccartoons}.

The periodic cartoon is then of the form
\[
f=\sum_{k \in \ZZ^2} 1_C(x+k), \quad C \subset \RR^2.
\]
Taking the Fourier transform, we observe that
\[
\hat{f} = \hat{1}_C \cdot \sum_{k \in \ZZ^2} e^{2 \pi i \langle k,\xi\rangle}
\]
with the second part being responsible for the sparsity in a Gabor system. This sparsity is the key for separation
from a single `large' cartoon. A further question for this particular situation might be: Given a set $C$ and a
lattice $A\ZZ^2$, how sparse is the corresponding periodization? However, this goes far beyond our model, wherefore
we do not treat this question here, but label it an interesting direction for future research.

\begin{figure}[ht]
\begin{center}
\vspace*{0.3cm}
\includegraphics[height=1.25in]{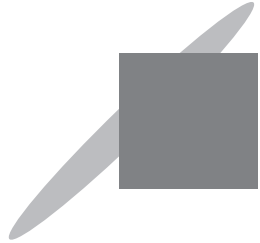}
\vspace*{-0.3cm}
\end{center}
\caption{The essential spatial support of a curvelet $\gamma_{j,k,\ell}$ and Gabor element $(g_s)_{m,n}$ satisfying the energy
balancing condition $s=s_j$ with $s_j$ defined in \eqref{eq:energybalancingexample}.}
\label{fig:criticalscale}
\end{figure}

Let us now take a closer look at Corollary \ref{coro:main}. To first build up intuition on the energy balancing condition,
we consider a curvelet of scale $j$ and a Gabor element of size $s$. If $s=s_j$ satisfies the energy balancing condition
\eqref{eq:energybalancingexample}, the spatial footprint of the Gabor element is about the size of the curvelet, since
$s_j = \omega(2^{j/2})$ and $s_j = o(2^j)$ in this case, which is illustrated in Figure \ref{fig:criticalscale}.

This also visualizes the relation between the filtered cartoon part $\cC_j$ and the filtered texture part $\cT_j (= \cT_{s_j,j})$,
since those will in fact be similarly connected. It is intuitively clear that this is the `worst scenario'. It seems
at first sight quite astonishing that separation is still possible in this case. However, the geometrical clustering of the
significant coefficients solves the mystery, and the notion of cluster coherence makes this precise. In fact, our result
shows that the whole cluster of Gabor elements $(g_{s_j})_{m,n}$ with $(m,n) \in B(0,r_{1,j}) \times B(0,r_{2,j})$
cannot be `glued together' so that it generates a single curvelet. The same is true for the whole cluster of significant
curvelets, which are unable to generate a single Gabor element. To visualize this would require to show the footprints
of curvelets and Gabor elements in a 5 dimensional `phase space': Spatial domain $\RR^2$, orientation $[0,2\pi]$, and
frequency $\RR^2$.


\subsection{Extensions}

Theorem \ref{theo:main} is amenable to the following generalizations and extensions.
\bitem
\item {\it Thresholding.}
Thresholding is a separation strategy, which is commonly utilized as a substitute for $\ell_1$ minimization
due to the fact that its algorithmic complexity is much lower. Hence one might ask: Can we derive a similar
asymptotic separation result when using thresholding? And in fact, this is true. We though decided not to include
the whole analysis, since it would go beyond the scope of this paper. We would just like to mention that
similar techniques as in \cite{DK08z} are employable, but the technical details, which will now in addition require
the special treatment of the scale of curvelets and the window size of Gabor atoms, will be very tedious.
\item {\it Intersections.}
Natural images are typically composed not only of one smoothly filled $C^2$ edge curve, but several with presumably
various intersections. In a similar way as it is commented upon in \cite{CD04}, we can also argue here that intersections
will not affect our analysis. Thus the framework can be extended to separating several cartoon-like objects
from texture.
\item {\it Other Systems.}
Theorem \ref{theo:main} holds without change for many other pairs of frames and bases, such as, e.g., by \cite{DK08b},
for the pair or shearlets (cf. \cite{GKL06,KL07,KKL10,KL10}) and Gabor systems.
\item {\it Noise.}
Theorem \ref{theo:main} is resilient to noise impact; an image composed of $\cC$ and $\cT$ with additive `sufficiently small' noise
exhibits the same asymptotic separation. This can be easily deduced by utilizing ideas from \cite[Sect. 8]{DK08a}.
\eitem

\section{Relative Sparsity}

\subsection{Cartoon}

Let $\cC_j$ be the filtered version of the cartoon-part $\cC$ of the considered
image, where -- as detailed in Subsection \ref{subsec:systems} -- we use $\hat{F}_j(\xi) = W(a_j|\xi|)$
as a filter. We now describe how to partition $\cC_j$ into smaller pieces which will
then in a second step be bent to `line objects' separately. The first part is inspired
by a technique employed in \cite{CD04}, whereas the second part follows ideas
introduced in \cite{DK08a}.

\begin{figure}
\centering
\includegraphics[height=1.75in]{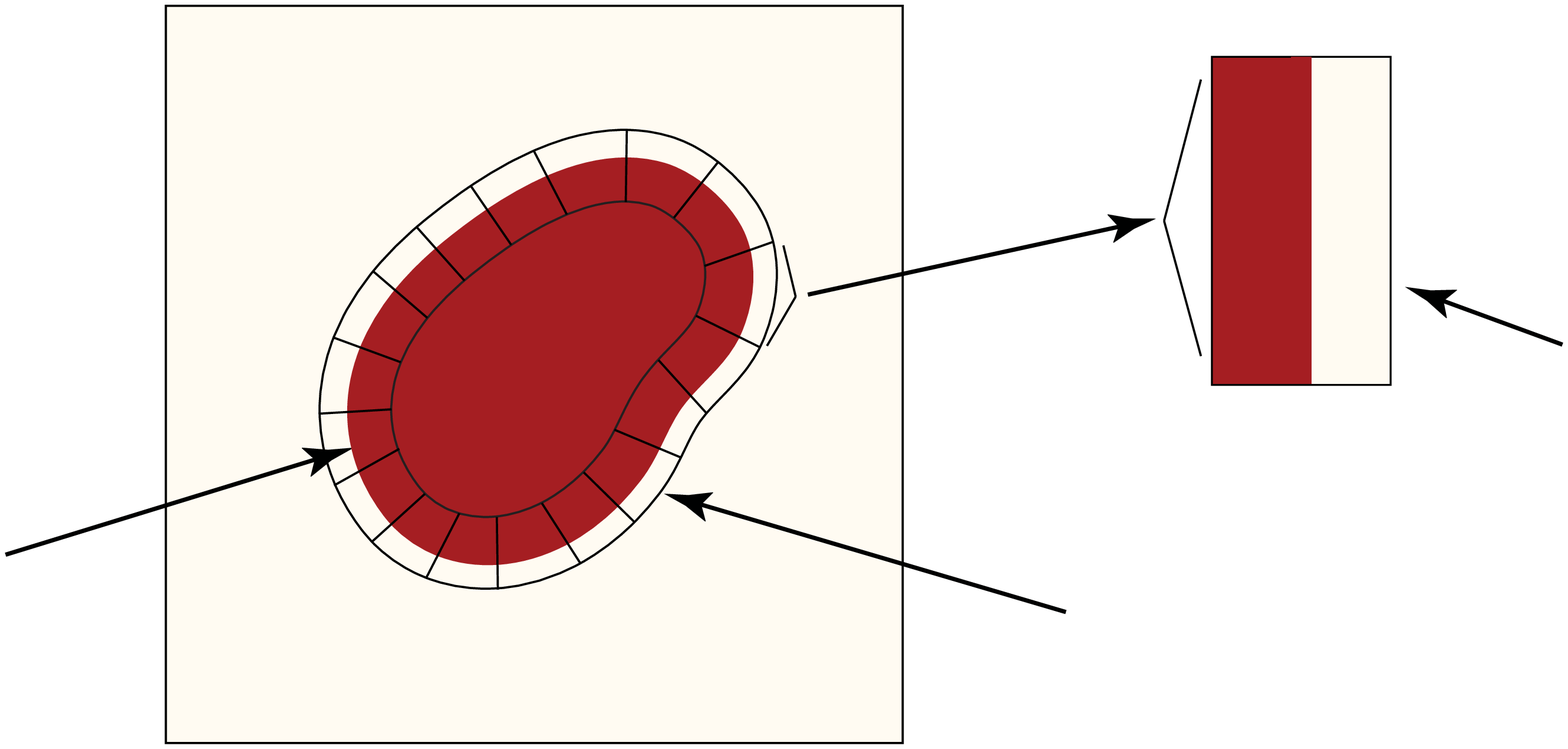}
\put(-316,29){$image(\tau)$}
\put(-83,18){$\cC_{j,Q}$}
\put(-108,87){$\phi_Q$}
\put(1,63){$H_{j,Q}$}
\caption{The filtered cartoon $\cC_j$ and the mapping $\phi_{Q} : \cC_{j,Q} \to H_{j,Q}$.}
\label{fig:tube}
\end{figure}

We start by first smoothly localizing $\cC_j$ near dyadic squares with a prescribed radius, which
is chosen such that the curvature of each piece is controllable, but independent on $j$.
More precisely, we define a partition of unity $(w_Q)_{Q \in \cQ}$ with those properties,
and let
\[
\cC_{j,Q} = \cC_j \cdot w_Q.
\]
We would like the reader to notice that the number of pieces $\cC_{j,Q}$ is finite and
independent on $j$. Also there is no need to consider those $\cC_{j,Q}$ which are smooth,
since their associated curvelet coefficients have sufficient decay. For each of the
remaining functions $\cC_{j,Q}$, without loss of
generality we can assume -- as was done in \cite[Sect. 6.1]{CD04} -- that the edge discontinuity
is centered at 0 and that its first derivative equals zero pointing in the vertical
direction. Hence, WLOG, the splitted cartoon might take a form as illustrated on the
LHS of Figure \ref{fig:tube}.

Having partitioned the filtered cartoon $\cC_j$, we now apply a diffeomorphism $\phi_{Q}$
to each piece $\cC_{j,Q}$, which equals the identity outside of a compact set, thereby
straightening out the discontinuity. For a similar strategy, we would like to refer the reader to
\cite[Sect. 7]{DK08a}. The resulting piece can now be modeled as
\[
H_{j,Q}=(H \star F_j) \cdot v_{Q},
\]
where $H = 1_{x_1 \ge 0}$ is the heaviside function and $v_{Q}$ is a $C^\infty$-function
supported in a dyadic square of sidelength $\rho$, say. An illustration of this slitting and
bending is presented in Figure \ref{fig:tube}.

We then require the following lemma, which is stated slightly more general than
needed. For this, we let $H^\nu$ be defined by
\[
\ip{\hat{H}^\nu}{f} = \int |\xi_1|^{-\nu} f(\xi_1,0) d\xi,
\]
hence, in particular, $H^1 = H$.

\begin{lemma}
\label{lem:HnuFj}
For each $N=1,2,\dots$, there is a constant $c_N$ so that
\[
  |\langle (H^\nu \star F_j) \cdot v_{Q}, \gamma_{a,b,\theta}  \rangle |
  \leq c_N \cdot a^{\nu - 1/4} \cdot  1_{|\log_2(a/a_j)|<3} \cdot \langle \frac{|\sin(\theta)-1|}{a} \rangle^{-N}
  \cdot \langle |b_1/a| \rangle^{-N}.
\]
\end{lemma}

\noindent
{\bf Proof.}
First, we observe that
\beq \label{eq:HnuFj1}
|\ip{(H^\nu \star F_j) \cdot v_{Q}}{ \gamma_{a,b,\theta}}| \le c_a 1_{|\log_2(a/a_j)|<3}.
\eeq
Therefore, we now restrict to the case $a = a_j$. Next, using Parseval, and letting
$\omega(\xi)$ denote the angular component of $\xi$,
\begin{eqnarray}\nonumber
\lefteqn{\ip{(H^\nu \star F_j) \cdot v_{Q}}{ \gamma_{a,b,\theta}}}\\ \nonumber
& = & \ip{\hat{H}^\nu \hat{F}_j}{\hat{v}_{Q} \star \hat{\gamma}_{a,b,\theta}}\\ \nonumber
& = & c \cdot \int |\xi_1|^{-\nu} \int \hat{v}_{Q}(\tau) a^{3/4} W(a|\xi_1-\tau|) V((\omega(\xi_1-\tau)-\theta)/\sqrt{a})
e^{-i\ip{b}{\xi_1-\tau}} d\tau\\ \label{eq:HnuFj2}
& & \cdot W(a|\xi_1|) d\xi_1.
\end{eqnarray}
By definition, the function $\hat{v}_{Q}$ is essentially supported in a square of radius $\rho^{-1}$.
Hence we approximately have $\xi_1-\tau \subseteq \{\xi : |\xi_2| \le \rho^{-1}\}$.
Assuming this condition, if $\theta \neq \pi/2$, then, for sufficiently small $a$,
\[
W(a|\xi_1-\tau|) V((\omega(\xi_1-\tau)-\theta)/\sqrt{a}) = 0 \qquad \mbox{for all } \xi_1, \tau.
\]
Thus, we can conclude that
\beq \label{eq:HnuFj3}
|\ip{(H^\nu \star F_j) \cdot v_{Q}}{ \gamma_{a,b,\theta}}| \le \langle \frac{|\sin(\theta)-1|}{a} \rangle^{-N}.
\eeq
Now, continuing the computation in \eqref{eq:HnuFj2} by using a change of variables, we obtain
\begin{eqnarray*}
\lefteqn{\ip{(H^\nu \star F_j) \cdot v_{Q}}{ \gamma_{a,b,\theta}}}\\
& = & c \cdot a^{\nu - 1/4} \cdot \int |\xi_1|^{-\nu} \int \hat{v}_{Q}(\tau)
W(|\xi_1-a\tau|) V((\omega(a^{-1}\xi_1-\tau)-\theta)/\sqrt{a}) e^{i\ip{b}{\tau}} d\tau\\
& & \cdot W(|\xi_1|) e^{-i (b_1/a)\xi_1} d\xi_1.
\end{eqnarray*}
Since $W$ is compactly supported, the size of the support of the function $\hat{G}_a$, defined by
\[
\hat{G}_a(\xi_1)=a^{\nu - 1/4} \cdot |\xi_1|^{-\nu} \int \hat{v}_{Q}(\tau)
W(|\xi_1-a\tau|) V((\omega(a^{-1}\xi_1-\tau)-\theta)/\sqrt{a}) e^{i\ip{b}{\tau}} d\tau W(|\xi_1|),
\]
is independent of $a$. Using repeated partial integration, we conclude
that
\[
|G(b)| \le \norm{\hat{G}_a}_\infty \cdot \langle |b_1/a| \rangle^{-N} \le a^{\nu - 1/4} \cdot \langle |b_1/a| \rangle^{-N}.
\]
The result is then completed by invoking \eqref{eq:HnuFj1} and \eqref{eq:HnuFj3}.
\qed

Let us now introduce the cluster of significant curvelet coefficients for $\cC_{j}$. First, we define
clusters for the manipulated pieces $H_{j,Q}$, which are then `lifted' to a cluster for $\cC_{j}$.
Intuitively, the cluster around $H_{j,Q}$ should spatially contain a carefully growing neighborhood
of the discontinuity $\{0\} \times [-\rho,\rho]$ and directionally contain again a carefully growing
neighborhood now of the direction of the discontinuity. For this, measuring the spatial distances
between a point $x$ and a set $A$ with
\[
d_2(x,A) = \min_{a \in A}\norm{x-a}_2, \quad x \in \bR^2,\,A \subseteq \bR^2,
\]
we define the neighborhood of the discontinuity $\{0\} \times [-\rho,\rho]$ in phase space by
\[
\cN_{2}^{PS}(a_j) = \{b \in \bR^2 : d_2(b,\{0\} \times [-2\rho,2\rho]) \le a_j^{1-\eps}\} \times [0,\sqrt{a_j}],
\]
where $\eps \in (0,\frac{1}{8})$. We are then led to define its corresponding neighborhood for discrete parameter
sets -- setting $\theta_{j,\ell} = \pi \ell / 2^{j/2}$ and $b_{j,k,\ell} = R_{\theta_{j,\ell}}D_{2^{-j}}k$ --
by
\[
 \t{\cS}_{j} = \{ (j,k,\ell) : (b_{j,k,\ell},\theta_{j,\ell}) \in \cN_{2}^{PS}(a_j)\}.
\]

In order to `lift' these clusters to a cluster for $\cC_{j}$, we will make use of the filtering matrix
associated with the filter $F_j$, i.e., of the matrix $M_{F_j} =
(\ip{\gamma_{\eta}}{F_j \star \gamma_{\eta'}})_{\eta,\eta'}$. We next need to recall that the action of
a diffeomorphism $\phi_Q$ on a distribution $f$ by $\phi_Q^\star f = f \circ \phi_Q$
induces a linear transformation on the space of curvelet coefficients. With $\alpha(f)$ the curvelet coefficients
of $f$ and $\beta(f)$ the curvelet coefficients of $\phi_Q^\star f$, we obtain
a linear operator $M_{\phi_Q}$ defined by
\[
    M_{\phi_Q} ( \alpha(f)) = \beta(f) .
\]
Letting now $t_{\eta',n_j^2}$ denote the amplitude of the $n_j^2$'th largest element of the $\eta'$'th column
of the matrix $M_{F_j} \cdot M_{(\phi_Q)^{-1}}$, we can define the overall cluster set for $\cC_{j}$ by
\[
   \cS_{1,j} = \bigcup_Q  \cS_{1,j,Q},
\]
where
\[
\cS_{1,j,Q} = \{ \eta : \eta' \in \t{\cS}_{j} \mbox{ and }
  |M_{F_j} \cdot M_{(\phi_Q)^{-1}}(\eta,\eta')| > t_{\eta',n_j^2} \}  .
\]

The following result will be needed for estimating the relative sparsity of
each $\cC_{j}$ with respect to the just defined cluster $\cS_{1,j}$.

\begin{lemma} \label{lem:analogLemma62}
Let $\alpha_j =  (\langle H_{j,Q} ,  \gamma_\eta  \rangle )_\eta$
denote the curvelet frame coefficients of $H_{j,Q}=(H \star F_j) \cdot v_{Q}$.
Then for all sufficiently large $j$,
\[
    \| \alpha_j \|_p \leq c_p  \cdot 2^{j(1/(2p)-3/4)}, \qquad \forall p > 0 .
\]
\end{lemma}

\noindent
{\bf Proof.}
This proof follows the lines of the proof of \cite[Lem. 6.2]{DK08a} very closely,
wherefore we decided not to state it explicitly. It should just be mentioned that
here we use Lemma \ref{lem:HnuFj} instead of \cite[Lem. 6.1]{DK08a}, which
leads to a difference in the asymptotic behavior by $2^{-j}$.
\qed

The decay estimate of the cluster approximate error $\delta_{1,j}$ is then
given by the following result.

\begin{lemma}\label{lem:delta2}
\[
 \delta_{1,j} = \sum_{\eta \in \Delta \setminus \cS_{1,j}} |\ip{\gamma_\eta}{\cC_j}|
 = o(2^{-j/2} ) , \qquad j \goto \infty.
\]
\end{lemma}

\noindent
{\bf Proof.}
Since the proof is very similar to \cite[Lem. 8.5]{DK08a}, we will not state it,
but merely mention that instead of \cite[Lem. 6.2]{DK08a} we here apply Lemma \ref{lem:analogLemma62}.
\qed


\subsection{Texture}
\label{sec:texture}

The decay estimate of the cluster approximation error $\delta_{2,j}$ is given by the following result.

\begin{lemma}\label{lem:delta1}
If there exist $r_{1,j}, r_{2,j} > 0$ such that
\beq\label{eq:condd1s}
\sum_{\stackrel{(m,n) \not\in B(0,r_{1,j}) \times B(0,r_{2,j})}{n \in \cA_{j}}} \sum_{\tilde{m} \in \ZZ^2} |d_{\tilde{m},n}| e^{-\frac{|\tilde{m}-m|}{2}}
= o(2^{-j/2}), \qquad j \to \infty,
\eeq
then
\[
 \delta_{2,j} = o(2^{-j/2}), \qquad j \to \infty.
\]
\end{lemma}

\noindent
{\bf Proof.}
First,
\[
\delta_{2,j} = \sum_{(m,n) \not\in B(0,r_{1,j})) \times B(0,r_{2,j}} |\ip{(g_{s_j})_{m,n}}{\cT_j}|.
\]
Similar to the proof of Lemma \ref{lem:estimateforTs}, we have
\[
|\ip{(g_{s_j})_{m,n}}{\cT_j}| \sim 1_{\cA_j}(n) \cdot \Big|\sum_{\tilde{m} \in \ZZ^2} d_{\tilde{m},n} e^{-\frac{|\tilde{m}-m|}{2}}\Big|, \qquad j \to \infty.
\]
This immediately implies the claim. \qed

To derive a better understanding, we analyze again the situation examined in Example \ref{exa:specialdmn}.

\begin{example} \label{exa:dmnford1j}
We here consider the situation that
\[
|d_{m,n}| \sim |m|^{-(2+\delta)} \cdot |n|^{-(2+\delta)}, \qquad \delta > 0,
\]
We aim to derive conditions on $r_{1,j}, r_{2,j} > 0$ such that \eqref{eq:condd1s} is satisfied, i.e., such that
\beq\label{eq:exaconddmn1}
T_1 \cdot T_2 + T_3 \cdot T_4 = o(2^{-j/2}) , \qquad j \to \infty,
\eeq
where
\begin{eqnarray*}
T_1 & = &  \sum_{|m| \ge r_{1,j}} \Big[C \cdot e^{-\frac{|m|}{2}} + \sum_{\tilde{m} \in \ZZ^2 \setminus \{0\}} |\tilde{m}|^{-(2+\delta)} e^{-\frac{|\tilde{m}-m|}{2}}\Big],\\
T_2 & = & \sum_{n \in \ZZ^2 \cap \cA_{j}} |n|^{-(2+\delta)},\\
T_3 & = & \sum_{|m| \le r_{1,j}} \Big[C \cdot e^{-\frac{|m|}{2}} + \sum_{\tilde{m} \in \ZZ^2 \setminus \{0\}} |\tilde{m}|^{-(2+\delta)} e^{-\frac{|\tilde{m}-m|}{2}}\Big],\\
T_4 & = & \sum_{n \in \ZZ^2 \cap \cA_{j} \cap B(0,r_{2,j})^c} |n|^{-(2+\delta)}.
\end{eqnarray*}

Firstly,
\[
C \cdot \int_{r_{1,j}}^\infty e^{-\frac{r}{2}}\cdot r dr + \int_{|x| \ge r_{1,j}} \int_{|y| \ge 1} |y|^{-(2+\delta)} e^{-\frac{|x-y|}{2}}
\sim r_{1,j}^{-\delta},
\]
hence
\beq \label{eq:exaconddmn2}
T_1 \sim r_{1,j}^{-\delta}.
\eeq
Secondly,
\[
\int_{\frac{2^{j-1}}{s_j}}^{\frac{2^{j+1}}{s_j}} r^{-(2+\delta)} \cdot r dr \sim \left(\frac{2^j}{s_j}\right)^{-\delta},
\]
hence
\beq \label{eq:exaconddmn3}
T_2 \sim \left(\frac{2^j}{s_j}\right)^{-\delta}.
\eeq
Thirdly,
\beq \label{eq:exaconddmn4}
T_3 \sim 1.
\eeq
And, fourthly,
\[
\int_{\max(\frac{2^{j-1}}{s_j},r_{2,j})}^{\frac{2^{j+1}}{s_j}} r^{-(2+\delta)} \cdot r dr
\sim \left\{\begin{array}{rcl}
0 & : & r_{2,j} = \omega\left(\frac{2^j}{s_j}\right),\\
\left(\frac{2^j}{s_j}\right)^{-\delta} & : & r_{2,j} = O\left(\frac{2^j}{s_j}\right).
\end{array} \right.
\]
hence
\beq \label{eq:exaconddmn5}
T_4
\sim \left\{\begin{array}{rcl}
0 & : & r_{2,j} = \omega\left(\frac{2^j}{s_j}\right),\\
\left(\frac{2^j}{s_j}\right)^{-\delta} & : & r_{2,j} = O\left(\frac{2^j}{s_j}\right).
\end{array} \right.
\eeq
Applying \eqref{eq:exaconddmn2}--\eqref{eq:exaconddmn5} to \eqref{eq:exaconddmn1} yields
\[
\delta_{1,j} \sim
r_{1,j}^{-\delta} \cdot \left(\frac{2^j}{s_j}\right)^{-\delta}
+\left\{\begin{array}{rcl}
0 & : & r_{2,j} = \omega\left(\frac{2^j}{s_j}\right),\\
\left(\frac{2^j}{s_j}\right)^{-\delta} & : & r_{2,j} = O\left(\frac{2^j}{s_j}\right).
\end{array} \right.
\]
Taking the energy matching condition \eqref{eq:energybalancingexample} from Example \ref{exa:specialdmn} into account,
\[
\delta_{1,j} \sim
r_{1,j}^{-\delta} \cdot 2^{-\frac{\delta}{2+2\delta}j}
+\left\{\begin{array}{rcl}
0 & : & r_{2,j} = \omega\left(\frac{2^j}{s_j}\right),\\
2^{-\frac{\delta}{2+2\delta}j} & : & r_{2,j} = O\left(\frac{2^j}{s_j}\right).
\end{array} \right.
\]
However, $2^{-\frac{\delta}{2+2\delta}j} = \omega(2^{-j/2})$. Hence, condition \eqref{eq:condd1s} is only fulfilled, if and
only if
\[
r_{1,j} = \omega\left(\frac{2^j}{s_j}\right)^{1/\delta} = \omega\left(2^{\frac{1}{2\delta(1+\delta)}j}\right)
\quad \mbox{and} \quad
r_{2,j} = \omega\left(\frac{2^j}{s_j}\right) = \omega\left(2^{\frac{1}{2(1+\delta)}j}\right).
\]
\end{example}


\section{Cluster Coherence}

\subsection{Interaction of Gabor Elements and Curvelets}

We start by estimating an inner product of a Gabor element and a curvelet. Depending on how $s_j$ relates to $2^j$, we
have to distinguish three cases. Note however, that the decay rate in all three cases is $O(2^{-j/4})$ as $j \to \infty$.

\begin{lemma} \label{lem:gamma_g_estimate}
For each $N=1,2,\dots$, there is  a constant $c_N$ so that,
\begin{enumerate}
\item[(i)] if $s_j = o(2^{j/2})$ as $j \to \infty$,
\[
  |\langle \gamma_{a_j,b,\theta} , (g_{s_j})_{m,n} \rangle | \leq c_N \cdot 2^{-3j/4} \cdot s_j \cdot 1_{M_{a_j,\theta},s_j}(n)  \cdot
    \langle |s_jb-\tfrac{m}{2}|_{R_\theta} \rangle^{-N},
\]
\item[(ii)] if $s_j = \Omega(2^{j/2})$ and $s_j = o(2^{j})$ as $j \to \infty$,
\[
  |\langle \gamma_{a_j,b,\theta} , (g_{s_j})_{m,n} \rangle | \leq c_N \cdot 2^{-j/4} \cdot 1_{M_{a_j,\theta},s_j}(n)  \cdot
    \langle |(s_jb_1-\tfrac{m_1}{2},2^{j/2}(b_2-\tfrac{m_2}{2s_j}))|_{R_\theta} \rangle^{-N},
\]
\item[(ii)] if $s_j = \Omega(2^{j/2})$ and $s_j = o(2^{j})$ as $j \to \infty$,
\[
  |\langle \gamma_{a_j,b,\theta} , (g_{s_j})_{m,n} \rangle | \leq c_N \cdot 2^{3j/4} \cdot s_j^{-1} \cdot 1_{M_{a_j,\theta},s_j}(n)  \cdot
    \langle |b-\tfrac{m}{2s_j}|_{D_{2^j}R_\theta} \rangle^{-N},
\]
\end{enumerate}
where $M_{a_j,\theta,s_j}$ is defined as in Theorem \ref{theo:main}.
\end{lemma}

\begin{figure}[ht]
\centering
\includegraphics[height=2.2in]{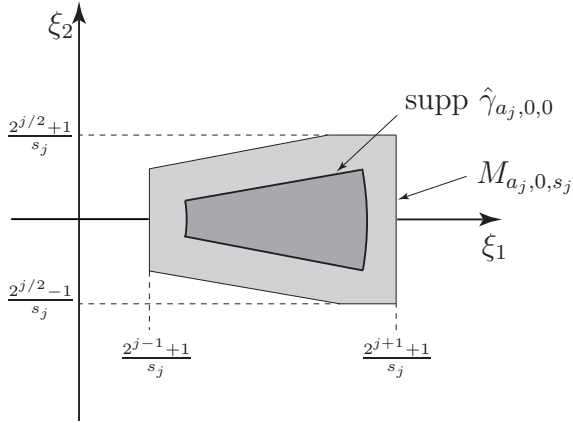}
\put(-10,90){\text{$M_{a_j,0,s_j}$}}
\put(-37,118){\text{supp $\hat{\gamma}_{a_j,0,0}$}}
\put(-54,22){\text{\tiny$\frac{2^{j+1}+1}{s_j}$}}
\put(-145,22){\text{\tiny$\frac{2^{j-1}+1}{s_j}$}}
\put(-188,106){\text{\tiny$\frac{2^{j/2}+1}{s_j}$}}
\put(-188,44){\text{\tiny$\frac{2^{j/2}-1}{s_j}$}}
\put(-8,63){$\xi_1$}
\put(-172,147){$\xi_2$}
\caption{Support of the set $M_{a_j,0,s_j}$.}
\label{fig:Matheta}
\end{figure}

\begin{figure}[ht]
\centering
\includegraphics[height=1.4in]{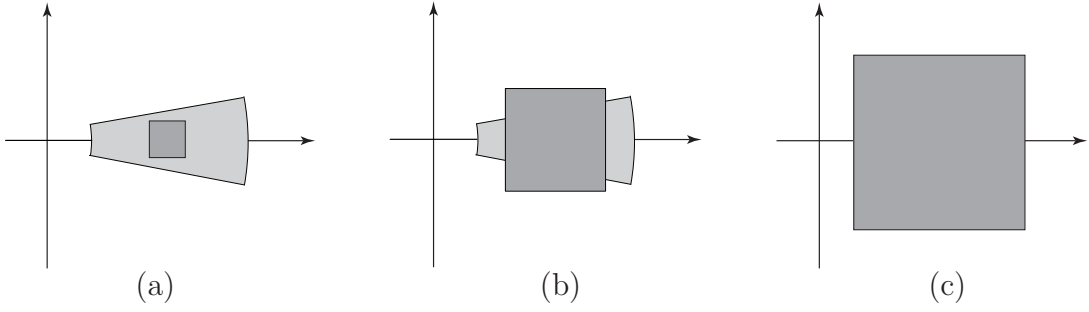}
\put(-360,-10){(a)}
\put(-207,-10){(b)}
\put(-60,-10){(c)}
\caption{Relation between the support of a curvelet $\hat{\gamma}_{j,0,k}$ and the support of a Gabor element $(\hat{g}_{s_j})_{m,(n_1,0)}$.
(a): Case $s_j = o(2^{j/2})$. (b): Case $s_j = \Omega(2^{j/2})$ and $s_j = o(2^{j})$. (c): Case $s_j = \Omega(2^{j})$.}
\label{fig:ThreeCases}
\end{figure}

Figure \ref{fig:Matheta} illustrates the set $M_{a_j,\theta,s_j}$ for $\theta = 0$, whereas the three cases in which
Lemma \ref{lem:gamma_g_estimate} is split into are illustrated in Figure \ref{fig:ThreeCases}.

\noindent
{\bf Proof of Lemma \ref{lem:gamma_g_estimate}.}
First, by Plancherel and the definition of $M_{a_j,\theta,s_j}$,
\begin{eqnarray}\nonumber
\lefteqn{\langle \gamma_{a_j,b,\theta},(g_s)_{m,n} \rangle}\\ \nonumber
& = & c \cdot a_j^{3/4} \int W(a_j|\xi|) V((\omega-\theta)/\sqrt{a_j}) \hat{g}_{s_j}(\xi-s_j n) e^{i(b-\frac{m}{2s_j})\xi} d\xi\\ \label{eq:coherence0}
& = & c \cdot a_j^{3/4} \cdot s_j^{-1} \cdot 1_{M_{a_j,\theta},s_j}(n) \int W(a_j|\xi|) V((\omega-\theta)/\sqrt{a_j}) \hat{g}(\tfrac{\xi}{s_j}-n) e^{i(b-\frac{m}{2s_j})\xi} d\xi.
\end{eqnarray}
It now remains to estimate
\[
T = \int W(a_j|\xi|) V((\omega-\theta)/\sqrt{a_j}) \hat{g}(\tfrac{\xi}{s_j}-n) e^{i(b-\frac{m}{2s_j})\xi} d\xi
\]
in the three cases the claim is split into. For this, in all three cases, WLOG we assume that $\theta=0$.

{\em Case $s_j = o(2^{j/2})$ as $j \to \infty$.}
With $\zeta = \frac{\xi}{s_j}$, and letting $\omega(\cdot)$ denote the angular component,
\[
|T| = \left|s_j^2 \cdot \int W(a_j s_j|\zeta|) V(\omega(s_j \zeta)/\sqrt{a_j}) \hat{g}(\zeta-n) e^{i(s_j b-\frac{m}{2})\zeta} d\zeta\right|.
\]
Applying integration by parts, for any $k=1, 2, ...$,
\begin{eqnarray*}
|T|
& = & s_j^2 \cdot |s_j b-\tfrac{m}{2}|^{-k} \left|\int \Delta^k[W(a_j s_j|\zeta|) V(\omega(s_j \zeta)/\sqrt{a_j}) \hat{g}(\zeta-n)] e^{i(s_j b-\frac{m}{2})\zeta} d\zeta\right|\\
& \le & s_j^2 \cdot |s_j b-\tfrac{m}{2}|^{-k} \int |\Delta^k[W(a_j s_j|\zeta|) V(\omega(s_j \zeta)/\sqrt{a_j}) \hat{g}(\zeta-n)]| d\zeta.
\end{eqnarray*}
Hence
\begin{eqnarray*}
(1+|s_j b-\tfrac{m}{2}|^{k}) \cdot |T|
& \le & \int |W(a_j s_j|\zeta|) V(\omega(s_j \zeta)/\sqrt{a_j}) \hat{g}(\zeta-n)| d\xi\\
& & + \int |\Delta^k[W(a_j s_j|\zeta|) V(\omega(s_j \zeta)/\sqrt{a_j}) \hat{g}(\zeta-n)]|d\zeta.
\end{eqnarray*}
It can be shown that, for each $k$, there exists $c_k < \infty$ such that, for all $a_j$,
\[
\int |W(a_j s_j|\zeta|) V(\omega(s_j \zeta)/\sqrt{a_j}) \hat{g}(\zeta-n)| d\xi + \hspace*{-0.1cm} \int |\Delta^k[W(a_j s_j|\zeta|) V(\omega(s_j \zeta)/\sqrt{a_j}) \hat{g}(\zeta-n)]| d\zeta
\le c_k.
\]
Further, for each $k =1, 2, ...$,
\[
\langle|s_j b-\tfrac{m}{2}|\rangle^k = (1+|s_j b-\tfrac{m}{2}|^2)^\frac{k}{2} \le \frac{k}{2}(1+|s_j b-\tfrac{m}{2}|^{k}),
\]
we obtain
\[
|T| \le c_N \cdot s_j^2 \cdot \langle|s_j b - \tfrac{m}{2} |\rangle^{-N}.
\]
Reinserting $\theta$ and using \eqref{eq:coherence0}, claim (i) is proved.

{\em Case $s_j = \Omega(2^{j/2})$ and $s_j = o(2^{j})$ as $j \to \infty$.}
With $\zeta = (s_j^{-1} \xi_1, 2^{-j/2} \xi_2) =: \varphi_{j}(\xi)$,
\begin{eqnarray*}
|T| & = & \left|2^{j/2} \cdot s_j \int W(a_j|\varphi_j^{-1}(\zeta)|) V(\omega(\varphi_j^{-1}(\zeta))/\sqrt{a_j}) \hat{g}(\zeta_1-n_1,2^{j/2}s_j^{-1}\zeta_2-n_2) \right.\\
& & \left.\cdot e^{i(s_jb_1-\tfrac{m_1}{2},2^{j/2}(b_2-\tfrac{m_2}{2s_j}))\zeta} d\zeta\right|.
\end{eqnarray*}
Repeated integration by parts -- similar as in the first case -- implies
\[
|T| \le c_N \cdot 2^{j/2}  \cdot s_j \cdot \langle|(s_jb_1-\tfrac{m_1}{2},2^{j/2}(b_2-\tfrac{m_2}{2s_j})) |\rangle^{-N}.
\]
Reinserting $\theta$ and using \eqref{eq:coherence0}, proves claim (ii).

{\em Case $s_j = \Omega(2^{j})$ as $j \to \infty$.}
With $\zeta = D_{2^{-j}}\xi$,
\[
|T| = \left|2^{3j/2} \int W(|\zeta|) V(\omega(D_{2^j} \zeta)/\sqrt{a_j}) \hat{g}(D_{2^j} \tfrac{\zeta}{s_j}-n) e^{i(b-\frac{m}{2s_j})D_{2^j}\zeta} d\zeta\right|.
\]
Repeated integration by parts -- similar as in the first case -- implies
\[
|T| \le c_N \cdot 2^{3j/2} \cdot \langle|b-\tfrac{m}{2s_j}|_{D_{2^j}}\rangle^{-N}.
\]
Claim (iii) now follows by reinserting $\theta$ and using \eqref{eq:coherence0}. \qed


\subsection{Cluster of Curvelets}

We can now prove negligible cluster coherence of $\cS_{1,j}$.

\begin{lemma} \label{lem:mu2}
\[
(\mu_c)_{1,j} = \mu_{c}(\cS_{1,j}, \{(g_{s_j})_\lambda\} ; \{ \gamma_\eta\}) \to 0,\quad \mbox{as } j \to \infty.
\]
\end{lemma}

\noindent
{\bf Proof.}
First, we observe that the definition of $\cS_{1,j}$ implies
\begin{eqnarray*}
(\mu_c)_{1,j}
& = & \max_{m,n} \sum_{(b,\theta) \in \cS_{1,j}} \absip{\gamma_{a_j,b,\theta}}{(g_{s_j})_{m,n}}\\
& = & \max_{m,n} \sum_Q \sum_{(b,\theta) \in \cS_{1,j,Q}} \absip{\gamma_{a_j,b,\theta}}{(g_{s_j})_{m,n}}\\
& \le & \#Q \cdot \max_Q \max_{m,n} \sum_{(b,\theta) \in \cS_{1,j,Q}} \absip{\gamma_{a_j,b,\theta}}{(g_{s_j})_{m,n}}
\end{eqnarray*}

Let us now first prove the claim in case that $\phi_Q = Id$, i.e., $M_Q = Id$, hence $\cS_{1,j,Q} = \t{\cS}_{j}$.
Notice that WLOG we can assume that $m=0$ and $n=(2^j,0) ( \in M_{a_j,\theta,s_j})$. Hence it remains to
continue the estimate
\beq \label{eq:Aimmuc}
(\mu_c)_{1,j} \le \#Q \cdot \max_Q \sum_{(b,\theta) \in \cS_{1,j,Q}} \absip{\gamma_{a_j,b,\theta}}{(g_{s_j})_{0,(2^j,0)}}.
\eeq
We now split the proof -- similar to the splitting in Lemma \ref{lem:gamma_g_estimate} -- in three cases.

{\em Case $s_j = o(2^{j/2})$ as $j \to \infty$.}
By \eqref{eq:Aimmuc} and Lemma \ref{lem:gamma_g_estimate},
\begin{eqnarray*}
(\mu_c)_{1,j} & \le & c_N \cdot 2^{-3j/4} \cdot s_j \sum_{(k,\ell) \in \t{\cS}_{j}} \langle |s_j(2^{-j} k_1, 2^{-j/2} k_2)|\rangle^{-N}\\
& \le & c_N \cdot 2^{-3j/4} \cdot s_j \cdot \pi^{-1} \cdot 2^{2 \eps j}  \sum_{|k_2| \le 2\rho + 2^{-j(1/2-\eps)}}
\langle |s_j 2^{-j/2} k_2|\rangle^{-N}\\
& \le & c_N \cdot 2^{-3j/4} \cdot s_j \cdot 2^{2 \eps j} \cdot s_j^{-1} \cdot (2\rho + 2^{-j(1/2-\eps)})\\
& \le & c_N \cdot 2^{j(2\eps - \frac34)}.
\end{eqnarray*}
Thus, since $\eps < \frac{1}{8}$,
\[
(\mu_c)_{1,j} \to 0,\quad \mbox{as } j \to \infty,
\]
which settles this case.

{\em Case $s_j = \Omega(2^{j/2})$ and $s_j = o(2^{j})$ as $j \to \infty$.}
By \eqref{eq:Aimmuc} and Lemma \ref{lem:gamma_g_estimate},
\begin{eqnarray*}
(\mu_c)_{1,j} & \le & c_N \cdot 2^{-j/4} \sum_{(k,\ell) \in \t{\cS}_{j}} \langle |(s_j b_1, 2^{j/2} b_2)|_{R_\theta}\rangle^{-N}\\
& \le & c_N \cdot 2^{-j/4} \cdot \pi^{-1} \cdot \sum_{|k_1| \le 2^{\eps j}}  \sum_{|k_2| \le 2\rho + 2^{-j(1/2-\eps)}}
\langle |(s_j 2^{-j} k_1, k_2)|\rangle^{-N}\\
& \le & c_N \cdot 2^{-j/4} \cdot 2^{2 \eps j}  \sum_{|k_2| \le 2\rho + 2^{-j(1/2-\eps)}}
\langle |k_2|\rangle^{-N}.
\end{eqnarray*}
Since the last sum is bounded by a constant,
\[
(\mu_c)_{1,j} \le c_N \cdot 2^{j(2\eps - \frac14)}.
\]
Since $\eps < \frac{1}{8}$,
\[
(\mu_c)_{1,j} \to 0,\quad \mbox{as } j \to \infty,
\]
which proves the claim for this case.

{\em Case $s_j = \Omega(2^{j})$ as $j \to \infty$.}
By \eqref{eq:Aimmuc} and Lemma \ref{lem:gamma_g_estimate},
\begin{eqnarray*}
(\mu_c)_{1,j} & \le & c_N \cdot 2^{3j/4} \cdot s_j^{-1} \sum_{(k,\ell) \in \t{\cS}_{j}} \langle |b|_{D_{2^j} R_\theta}\rangle^{-N}\\
& \le & c_N \cdot 2^{3j/4} \cdot s_j^{-1} \cdot \pi^{-1} \cdot \sum_{k} \langle |k|\rangle^{-N}\\
& \le & c_N \cdot 2^{3j/4} \cdot s_j^{-1}.
\end{eqnarray*}
Hence,
\[
(\mu_c)_{1,j} \to 0,\quad \mbox{as } j \to \infty.
\]
Hence also this case is proved.

The general claim now follows by using similar arguments as in \cite[Lem. 8.8]{DK08a}.
\qed


\subsection{Cluster of Gabor Elements}

Finally, we arrive at the study of the cluster coherence of $\cS_{2,j}$.

\begin{lemma} \label{lem:mu1}
Suppose that one of conditions (i) and (ii) is satisfied:
\begin{enumerate}
\item[(i)] Suppose $s_j = o(2^{j})$ as $j \to \infty$, and $r_{2,j}$ is chosen such that
\[
|B(0,r_{2,j}) \cap M_{a_j,0,s_j} \cap \ZZ^2| = o(2^{3j/4} \cdot s_j^{-1}),\quad \mbox{as } j \to \infty.
\]
\item[(ii)] Suppose $s_j = \Omega(2^{j})$ as $j \to \infty$, and $r_{1,j}$ is chosen such that
\[
|B(0,r_{1,j}) \cap \ZZ^2| = o(2^{-3j/4} \cdot s_j),\quad \mbox{as } j \to \infty.
\]
\end{enumerate}
Then
\[
(\mu_c)_{2,j} = \mu_{c}(\cS_{2,j},  \{ \gamma_\eta\} ; \{(g_{s_j})_\lambda\} ) \to 0,\quad \mbox{as } j \to \infty.
\]
\end{lemma}

\noindent
{\bf Proof.}
First,
\begin{eqnarray}\nonumber
(\mu_c)_{2,j} & = & \max_{b,\theta} \sum_{(m,n) \in \cS_{2,j}} \absip{\gamma_{a_j,b,\theta}}{(g_{s_j})_{m,n}}\\ \label{eq:coh1j1}
& = & \sum_{(m,n) \in \cS_{2,j}} \absip{\gamma_{a_j,0,0}}{(g_{s_j})_{m,n}},
\end{eqnarray}
since, for symmetry reasons, WLOG we can assume that the maximum is attained in $b=0$ and $\theta = 0$.
The two cases will now be dealt with separately. We start with the first, which -- similar to the splitting
in Lemma \ref{lem:gamma_g_estimate} -- we separate into two subcases.

{\em Case $s_j = o(2^{j/2})$ as $j \to \infty$.}
By \eqref{eq:coh1j1} and Lemma \ref{lem:gamma_g_estimate},
\[
(\mu_c)_{2,j} \le c_N \cdot 2^{-3j/4} \cdot s_j \sum_{m \in B(0,r_{1,j})} \langle |\tfrac{m}{2}|\rangle^{-N} \sum_{n \in B(0,r_{2,j}) \cap M_{a_j,0,s_j} \cap \ZZ^2} 1.
\]
Since the first sum is bounded by a constant and by condition (i),
\[
(\mu_c)_{2,j} = o(1),\quad \mbox{as } j \to \infty,
\]
this case is settled.

{\em Case $s_j = \Omega(2^{j/2})$ and $s_j = o(2^{j})$ as $j \to \infty$.}
By \eqref{eq:coh1j1} and Lemma \ref{lem:gamma_g_estimate},
\[
(\mu_c)_{2,j} \le c_N \cdot 2^{-j/4} \sum_{m \in B(0,r_{1,j})} \langle |(\tfrac{m_1}{2},2^{j/2} \tfrac{m_2}{2s_j}|\rangle^{-N}
\sum_{n \in B(0,r_{2,j}) \cap M_{a_j,0,s_j} \cap \ZZ^2} 1.
\]
Since
\[
\int_0^{r_{1,j}} \langle|2^{j/2} s_j^{-1} x|\rangle^{-N} dx = 2^{-j/2} \cdot s_j \int_0^{2^{j/2} s_j^{-1} r_{1,j}} \langle|x|\rangle^{-N} dx
\sim 2^{-j/2} \cdot s_j,
\]
we obtain -- also exploiting condition (i),
\[
(\mu_c)_{2,j} \le c_N \cdot 2^{-3j/4} \cdot s_j \sum_{n \in B(0,r_{2,j}) \cap M_{a_j,0,s_j} \cap \ZZ^2} 1
= o(1),\quad \mbox{as } j \to \infty.
\]
This proves the claim for this case.

{\em Case $s_j = \Omega(2^{j})$ as $j \to \infty$.}
By \eqref{eq:coh1j1} and Lemma \ref{lem:gamma_g_estimate},
\[
(\mu_c)_{2,j} \le c_N \cdot 2^{3j/4} \cdot s_j^{-1} \sum_{m \in B(0,r_{1,j})} \langle |b-\tfrac{m}{2s_j}|_{D_{2^j}}\rangle^{-N}
\sum_{n \in B(0,r_{2,j}) \cap M_{a_j,0,s_j} \cap \ZZ^2} 1.
\]
Since $| M_{a_j,0,s_j} \cap \ZZ^2| \le const$,
\begin{eqnarray*}
(\mu_c)_{2,j} & \le & c_N \cdot 2^{3j/4} \cdot s_j^{-1} \sum_{m \in B(0,r_{1,j})} \langle |\tfrac12 (2^j s_j^{-1} m_1, 2^{j/2} s_j^{-1} m_2)|\rangle^{-N}\\
& \le & c_N \cdot 2^{3j/4} \cdot s_j^{-1} |\{m \in \ZZ^2 : |m| \le r_{1,j}\}|.
\end{eqnarray*}
By condition (ii),
\[
(\mu_c)_{2,j} = o(1),\quad \mbox{as } j \to \infty,
\]
The lemma is proved. \qed


\section{Proofs}

\subsection{Proofs of Results from Section \ref{sec:intro}}
\label{subsec:proofs_1}

\subsubsection{Proof of Lemma \ref{lem:estimateforCj}}

\noindent
{\bf Proof.}
We first observe that, for $W$ sufficiently nice,
\beq \label{eq:estimatecj}
\| \cC_j\|_2^2 = \int |W(a_j|\xi|)|^2 |\hat{\cC}_j(\xi)|^2 d\xi \sim \int_{\tilde{\cA}_j} |\hat{\cC}(\xi)|^2 d\xi,
\eeq
where $\tilde{\cA}_j = \{\xi : 2^j \le |\xi| < 2^{j+1}\}$. Using \eqref{eq:extraconditiononC}, we can continue \eqref{eq:estimatecj}
to
\[
\| \cC_j\|_2^2 \sim \int_{2^j}^{2^{j+1}} r^{-2} dr \sim 2^{-j}. \quad \mbox{\qed}
\]


\subsubsection{Proof of Lemma \ref{lem:estimateforTs}}

First, by the change of variable $\omega = \xi/s$ and the support condition on $\hat{g}$,
\begin{eqnarray} \nonumber
\| \cT_{s,j}\|_2^2
& = & \sum_{m,n} \sum_{\tilde{m},\tilde{n}} \int d_{m,n} \overline{d_{\tilde{m},\tilde{n}}} W^2(a_j|\xi|) \hat{g}_s(\xi-s n)\overline{\hat{g}_s(\xi-s \tilde{n})}
e^{i \frac{m-\tilde{m}}{2s}\xi} d\xi\\ \nonumber
& = & s^2 \cdot \sum_{\stackrel{m,n,\tilde{m},\tilde{n}}{|n-\tilde{n}| \le 1}} \int d_{m,n} \overline{d_{\tilde{m},\tilde{n}}} W^2(a_js|\omega|) \hat{g}_s(s(\omega-n))
\overline{\hat{g}_s(s(\omega-\tilde{n}))} e^{i \frac{m-\tilde{m}}{2}\omega} d\omega\\ \label{eq:estimateTs1}
& = & \sum_{\stackrel{m,n,\tilde{m},\tilde{n}}{|n-\tilde{n}| \le 1}} d_{m,n} \overline{d_{\tilde{m},\tilde{n}}}  \int W^2(a_js|\omega|) \hat{g}(\omega-n)
\overline{\hat{g}(\omega-\tilde{n})} e^{i \frac{m-\tilde{m}}{2}\omega} d\omega.
\end{eqnarray}
For `sufficiently nice' $W$,
\[
\int W^2(a_js|\omega|) \hat{g}(\omega-n) \overline{\hat{g}(\omega-\tilde{n})} e^{i \frac{m-\tilde{m}}{2}\omega} d\omega
\sim \int_{\cA_{s,j}} \hat{g}(\omega-n) \overline{\hat{g}(\omega-\tilde{n})} e^{i \frac{m-\tilde{m}}{2}\omega} d\omega,
\]
hence, continuing \eqref{eq:estimateTs1}, and taking into account that for each $n$, there exist a finite number of
$\tilde{n}$'s independent on $j$ satisfying $|n-\tilde{n}| \le 1$,
\begin{eqnarray} \nonumber
\|\cT_{s,j}\|_2^2
& \sim & \sum_{\stackrel{m,n,\tilde{m},\tilde{n}}{|n-\tilde{n}| \le 1}} d_{m,n} \overline{d_{\tilde{m},\tilde{n}}}
\int_{\cA_{s,j}} \hat{g}(\omega-n) \overline{\hat{g}(\omega-\tilde{n})} e^{i \frac{m-\tilde{m}}{2}\omega} d\omega\\ \nonumber
& = & \sum_{m,\tilde{m}} \sum_{\stackrel{n,\tilde{n} \in \ZZ^2 \cap \cA_{s,j}}{|n-\tilde{n}| \le 1}} d_{m,n} \overline{d_{\tilde{m},\tilde{n}}}
\int_{\cA_{s,j}} \hat{g}(\omega-n) \overline{\hat{g}(\omega-\tilde{n})} e^{i \frac{m-\tilde{m}}{2}\omega} d\omega\\ \label{eq:estimateTs2}
& \sim & \sum_{m,\tilde{m}} \sum_{n \in \ZZ^2 \cap \cA_{s,j}} d_{m,n} \overline{d_{\tilde{m},n}}
\int_{\cA_{s,j}} |\hat{g}(\omega-n)|^2 e^{i \frac{m-\tilde{m}}{2}\omega} d\omega.
\end{eqnarray}
Now, WLOG we assume that the support of $\hat{g}$ is always entirely contained in $\cA_{s,j}$; then, by \eqref{eq:estimateTs2},
\beq \label{eq:estimateTs3}
\|\cT_{s,j}\|_2^2
\sim  \sum_{m,\tilde{m}} \sum_{n \in \ZZ^2 \cap \cA_{s,j}} d_{m,n} \overline{d_{\tilde{m},n}}
\int_{\RR^2} |\hat{g}(\omega-n)|^2 e^{i \frac{m-\tilde{m}}{2}\omega} d\omega.
\eeq
Since $|g(x)| \sim e^{-|x|}$,
\begin{eqnarray*}
\int_{\RR^2} |\hat{g}(\omega-n)|^2 e^{i \frac{m-\tilde{m}}{2}\omega} d\omega
& = & e^{i \frac{m-\tilde{m}}{2}n} \cdot \int_{\RR^2} |\hat{g}(\omega)|^2 e^{i \frac{m-\tilde{m}}{2}\omega} d\omega\\
& = & e^{i \frac{m-\tilde{m}}{2}n} \cdot (g \star g^\star)(\tfrac{m-\tilde{m}}{2})\\
& \sim & e^{-\frac{|m-\tilde{m}|}{2}}.
\end{eqnarray*}
Thus, by \eqref{eq:estimateTs3},
\[
\|\cT_{j}\|_2^2 \sim \sum_{m,\tilde{m}} \sum_{n \in \ZZ^2 \cap \cA_{s,j}} d_{m,n} \overline{d_{\tilde{m},n}} e^{-\frac{|m-\tilde{m}|}{2}},
\]
and the lemma is proved. \qed


\subsubsection{Proof of Theorem \ref{theo:main}}

The theorem follows immediately from applying Lemmata \ref{lem:delta1}, \ref{lem:delta2}, \ref{lem:mu1},
and \ref{lem:mu2} to Proposition \ref{prop:coherenceestimate}. \qed


\subsubsection{Proof of Corollary \ref{coro:main}}

We choose $r_{1,j}> 0$ such that
\[
r_{1,j} = \omega(2^{\frac{1}{2\delta(1+\delta)}j}), \quad \mbox{as } j \to \infty,
\]
and  $r_{2,j}> 0$ such that
\beq \label{eq:proofcoro0}
r_{2,j} = \omega(2^{\frac{1}{2(1+\delta)}j}), \quad \mbox{as } j \to \infty.
\eeq

It now suffices to check the sufficient conditions posed in Theorem \ref{theo:main}. Firstly,
by Example \ref{exa:dmnford1j},
\[
\sum_{\stackrel{(m,n) \not\in B(0,r_{1,j}) \times B(0,r_{2,j})}{n \in \cA_{j}}} \sum_{\tilde{m} \in \ZZ^2} |d_{\tilde{m},n}| e^{-\frac{|\tilde{m}-m|}{2}}
= o(2^{-j/2}), \qquad j \to \infty.
\]
Now notice that the energy matching condition in this case, i.e., \eqref{eq:energybalancingexample}, implies that
$s_j = \Omega(2^{j/2})$ and $s_j = o(2^{j})$ as $j \to \infty$. Hence we need to prove
\beq \label{eq:proofcoro1}
|B(0,r_{2,j}) \cap M_{a_j,0,s_j} \cap \ZZ^2| = o(2^{3j/4} \cdot s_j^{-1}) = o(2^{\frac{1-\delta}{4(1+\delta)}j}),\quad \mbox{as } j \to \infty.
\eeq
For this, by \eqref{eq:proofcoro0},
\[
|B(0,r_{2,j}) \cap M_{a_j,0,s_j} \cap \ZZ^2| = |M_{a_j,0,s_j} \cap \ZZ^2|
\sim s_j^{-1} \cdot 2^j \cdot s_j^{-1} \cdot 2^{j/2} = s_j^{-2} \cdot 2^{3j/2},\quad \mbox{as } j \to \infty.
\]
Since, by \eqref{eq:energybalancingexample},
\[
s_j^{-2} \cdot 2^{3j/2} = 2^{\frac{1-\delta}{2(1+\delta)}j},
\]
claim \eqref{eq:proofcoro1} is satisfied if and only if $\delta > 1$, which was assumed.

The claim now follows from Theorem \ref{theo:main}. \qed


\end{document}